\documentclass{article} 
\usepackage[usenames,dvipsnames]{color}
\usepackage{amsmath} 
\usepackage{amssymb} 
\usepackage{bm}
\usepackage{amsthm}
\usepackage{CJK}
\usepackage{indentfirst}
\usepackage{tikz}
\usepackage{amsmath,amsfonts,amsthm,mathrsfs}
\usepackage{hyperref}
\numberwithin{equation}{section}
\usepackage{geometry}
\usepackage{extarrows}
\usepackage{setspace}

\numberwithin{equation}{section}
\newtheorem{theorem}{Theorem}[section]
\newtheorem{lemma}{Lemma}[section]

\newtheorem{definition}{Definition}[section]

\usepackage{algorithm,algpseudocode}
\usepackage{amsthm}
\usepackage{esint} 
\usepackage{subfigure}
\usepackage{natbib}
\usepackage{lineno}

\def\NN{\mathbb{N}}
\def\PP{\mathbb{P}}
\def\RR{\mathbb{R}}

\def\SS{\mathbb{S}}

\def\mal{\max\limits}
\def\mil{\min\limits}
\def\sul{\sum\limits}

\makeatletter
\renewcommand\@fnsymbol[1]{\arabic{footnote}}
\makeatother

\title{Parameter Bounds and Inverse Approximation for\\ Linearized ReLU$^k$ Networks on the Sphere}

\author{Yuqing Liu\thanks{RPTU University Kaiserslautern-Landau, Kaiserslautern and Landau in der Pfalz, Rhineland-Palatinate, Germany}\and Tong Mao\thanks{King Abdullah University of Science and Technology, Thuwal 23955, Saudi Arabia}\and Jonathan Siegel\thanks{Texas A\&M University, College Station, TX 77843, USA}\and Jinchao Xu\thanks{King Abdullah University of Science and Technology, Thuwal 23955, Saudi Arabia}}
\date{}

\begin{document}
	
\maketitle

\begin{abstract}
We study direct and inverse approximation for linearized shallow ReLU$^k$ networks on $\SS^d$ with quasi-uniform hidden parameters. For $0<\alpha<\frac{d+2k+1}{2d}$, under the coefficient growth condition $\bigl(n\sum_{j=1}^na_{j,n}^2\bigr)^{\frac{1}{2}}\leq n^\alpha M$, we prove that approximation at the order $n^{-\frac{d+2k+1}{2d}+\alpha}$ is equivalent to membership in the parity-compatible Besov space $\mathcal B_{2,\infty,*}^{\frac{d+2k+1}{2}-\alpha d}(\SS^d)$. Mere convergence under the same coefficient budget does not imply this regularity. For unconstrained approximation along a nested parameter sequence, a single ridge function gives a counterexample when $\alpha<\frac{1}{2}$, whereas for $\alpha>\frac{1}{2}$ the endpoint Besov inverse conclusion follows from a Bernstein inequality under antipodal quasi-uniformity. The inverse conclusion is also recovered when the prescribed rate holds for two quasi-uniform parameter sequences whose union is antipodally separated.
\end{abstract}

\section{Introduction}
The approximation theory of shallow neural networks seeks to quantify how the approximation error depends on the network width, the activation function, the regularity of the target, and the size of the representing coefficients. For Barron and variation-type spaces, dimension-independent rates are commonly obtained through probabilistic or greedy discretization arguments \cite{Jones1992,barron1993universal,barron1994approximation,devore1996some,bach2017breaking,klusowski2018approximation,siegel2022high,siegel2022sharp,siegel2022optimal,mhaskar2023tractability}. For Sobolev classes, the rate instead reflects both smoothness and dimension \cite{mhaskar1993approximation,petrushev1998approximation,pinkus1999approximation}. These direct estimates give sufficient conditions for attaining a prescribed rate. The converse question---what regularity is forced by that rate---is subtler and is closely related to Bernstein inequalities, inverse approximation theorems, and the structure of neural-network approximation classes \cite{siegel2022sharp,siegel2023characterization,guo2020optimal,E2019,bernstein2020distance,huang2023generalization}.

We focus on the ReLU$^k$ activation $\sigma_k(t)=(\max\{0,t\})^k$, $k\in\NN_0$, with the usual Heaviside convention when $k=0$, and on linearized shallow networks, in which the hidden parameters are fixed and only the output coefficients vary. By the positive homogeneity of $\sigma_k$, the parameter variables may be normalized to directions on the unit sphere. We consider the approximation problem on
\[
    \SS^d=\{\eta\in\RR^{d+1}:|\eta|=1\},
\]
where the basic ridge functions take the form $\sigma_k(\theta\cdot\eta)$ with $\theta,\eta\in\SS^d$.

The spherical setting is useful both as a normalized model for affine ReLU$^k$ networks and as an approximation problem in its own right. Spherical neural networks arise naturally for data with directional or rotational structure \cite{cohen2018spherical,jiang2019spherical,defferrard2019deepsphere}, and their approximation and generalization properties have been studied in \cite{fang2020theory,feng2023generalization,yang2026spherical}. Harmonic analysis on the sphere also provides localized polynomial kernels, cubature formulas, and stable discretizations on quasi-uniform point sets \cite{mhaskar2001spherical,mhaskar2004polynomialLocal,mhaskar2005sphereLocalized,mhaskar2006weighted,mhaskar2019reluZonal,mhaskar2020kernel}. These tools connect the geometry of the hidden parameters with approximation accuracy and coefficient stability.

Let $\Theta_n=\{\theta_{j,n}^*\}_{j=1}^n\subset\SS^d$. For $\alpha\geq0$ and $M>0$, define
\begin{equation}\label{eqn_def_LnM}
    L_{n,M}^{k,\alpha}\bigl(\Theta_n\bigr)=\Big\{g\in\mathcal L^2(\SS^d):~g(\eta)=\sum_{j=1}^na_{j,n}\sigma_k\bigl(\theta_{j,n}^*\cdot\eta\bigr),\quad\Bigl(n\sul_{j=1}^na_{j,n}^2\Bigr)^{\frac{1}{2}}\leq n^\alpha M\Big\}.
\end{equation}
We write $L_n^k(\Theta_n)$ when no constraint is imposed on the coefficients. The constant $M$ is always independent of $n$; in a direct approximation statement for a given $f$, it may depend on $f$. Schwartz inequality shows the coefficient bound of the Barron space \cite{siegel2023characterization} is consistent with the constraint in \eqref{eqn_def_LnM} when $\alpha=0$:
\begin{equation*}
    \sul_{j=1}^n|a_{j,n}|\leq\Bigl(n\sul_{j=1}^na_{j,n}^2\Bigr)^{\frac{1}{2}}\leq M
\end{equation*}

The identity
\[
    \sigma_k(\theta\cdot\eta)+(-1)^k\sigma_k(-\theta\cdot\eta)=(\theta\cdot\eta)^k
\]
shows that every element of $L_n^k(\Theta_n)$ has a prescribed parity up to a spherical polynomial of degree at most $k$. Therefore, for a function space $\mathcal S(\SS^d)$ containing $\PP_k(\SS^d)$, we write
\begin{equation*}
    \mathcal S_*(\SS^d):=\Big\{f\in\mathcal S(\SS^d):~\Pi_mf=0\ \text{for every }m>k\text{ with }m-k\text{ even}\Big\}.
\end{equation*}
Here $\Pi_m$ denotes the spherical harmonic projection introduced below. Equivalently, up to $\PP_k(\SS^d)$, $\mathcal S_*(\SS^d)$ consists of even functions when $k$ is odd and odd functions when $k$ is even.

A sharp direct approximation result was established in \cite{liu2025integral}. For a quasi-uniform parameter sequence $\{\Theta_n\}_{n=1}^\infty$ and $0\leq \alpha\leq\frac{d+2k+1}{2d}$, one has, with $M$ controlled by the Sobolev norm of $f$,
\begin{equation}\label{eqn_appr_rate_main}
    \inf_{f_n\in L_{n,M}^{k,\alpha}(\Theta_n)}\|f-f_n\|_{\mathcal L^2(\SS^d)}\lesssim n^{-\frac{d+2k+1}{2d}+\alpha}\|f\|_{\mathcal H_*^{\frac{d+2k+1}{2}-d\alpha}(\SS^d)}.
\end{equation}
The endpoint $\frac{d+2k+1}{2}$ is determined by the decay of the spherical harmonic coefficients of a ReLU$^k$ ridge function.

Two converse results provide the immediate background for this work. First, the integral characterization in \cite{liu2025integral} shows that if, for some quasi-uniform sequence $\{\Theta_n\}_{n=1}^\infty$,
\begin{equation}\label{eqn:appr_equiv_sob}
    \lim_{n\to\infty}\inf_{f_n\in L_{n,M}^{k,0}(\Theta_n)}\|f-f_n\|_{\mathcal L^2(\SS^d)}=0,
\end{equation}
then $f\in\mathcal H_*^{\frac{d+2k+1}{2}}(\SS^d)$. Second, the Bernstein inequality established in \cite{liu2026optimal} gives an inverse theorem for the unconstrained spaces. If $\{\Theta_n\}_{n\geq1}$ is nested and antipodally quasi-uniform and $\frac{1}{2}<\alpha<\frac{d+2k+1}{2d}$, then
\begin{equation}\label{eqn_decay_inv_appr}
    \inf_{f_n\in L_n^k(\Theta_n)}\|f-f_n\|_{\mathcal L^2(\SS^d)}\lesssim n^{-\frac{1}{2}-\frac{2k+1}{2d}+\alpha},
\end{equation}
implies $f\in\mathcal H_*^r(\SS^d)$ for every $0\leq r<\frac{d+2k+1}{2}-\alpha d$. We strengthen these subcritical Sobolev conclusions below to the endpoint Besov space.

Direct and inverse estimates of this kind are widespread in classical
approximation theory: together they determine whether a rate is sharp and
whether improved convergence forces additional regularity of the target
\cite{bernstein1912quelques,nikolskii1975approximation,
devore1989optimal,zhou2006approximation,maiorov2006pseudo}. One of the key concepts in this theory is the notion of an approximation space (see \cite{devore1998nonlinear}) corresponding to a sequence of linear spaces. Specifically, let $V_n\subset \mathcal L^2(\SS^d)$ be a linear space of dimension $n$. For a given decay rate $s > 0$, the approximation spaces with respect to the sequence of linear spaces $V_n$ is defined by the semi-norm
\begin{equation}
    |f|_{\mathcal{A}^s(\{V_n\}_{n=1}^\infty)} := \sup_{n \geq 1}n^\alpha E(f,V_n),
\end{equation}
where
\begin{equation}
    E(f,V_n) = \inf_{f_n\in V_n} \|f - f_n\|_{\mathcal L^2(\SS^d)}
\end{equation}
is the error of best approximation by the space $V_n$.

A useful classical comparison to our results is provided by the setting of linear approximation by splines. The simplest case of this is linear approximation by piecewise constants, which is detailed in Section 3.1 in \cite{devore1998nonlinear}, for instance. In this theory, the approximation spaces for a sequence of (non-nested) piecewise constant linear spaces $V_n$ are shown to be Besov spaces, but only for decay rates up to $n^{-s}$ with $s < 1$ (see \cite{devore1993constructive}, Chapter 7, Section 7). This is due to the fact that higher-order smoothness cannot be captured by piecewise constants, and indeed for higher order decay rates the approximation spaces cannot be described in a simple form.

When we consider the linear spaces $V_n := L_n^k(\Theta_n)$ for a sequence of quasi-uniform points $\Theta_n$, one might expect that the same characterization of the approximation spaces holds. Our results show that this is indeed the case, up to approximation order $s < \frac12+\frac{2k+1}{2d}$, provided that we additionally constrain the coefficients of the approximant as in the definition \eqref{eqn_def_LnM}, i.e., we consider
\begin{equation}
    E(f,L_{n,M}^{k,\alpha}(\Theta_n)) := \inf_{f_n\in L_{n,M}^{k,\alpha}(\Theta_n)} \|f - f_n\|_{\mathcal L^2(\SS^d)},
\end{equation}
instead of $E(f,V_n)$ for the linear spaces $V_n = L_n^k(\Theta_n)$.
This additional constraint is required since the ReLU$^k$ ridge functions are globally supported and highly correlated, in contrast to the spline basis. Motivated by this comparison, we address the following two questions:

\begin{enumerate}
    \item Which function classes can be approximated when a prescribed growth bound is imposed on the coefficients?
    \item Which regularity properties can be inferred from a prescribed approximation rate, with or without an explicit coefficient bound?
\end{enumerate}
For $\alpha\geq0$, the natural coefficient and approximation scales considered here are
\begin{equation}
    \Bigl(n\sul_{j=1}^na_{j,n}^2\Bigr)^{\frac{1}{2}}\leq n^\alpha M,\qquad\Big\|f-\sul_{j=1}^na_{j,n}\sigma_k(\theta_{j,n}^*\cdot\circ)\Big\|_{\mathcal L^2(\SS^d)}\lesssim n^{-\frac{1}{2}-\frac{2k+1}{2d}+\alpha}.
\end{equation}
For $\alpha=0$, these are the coefficient scale and endpoint rate in \eqref{eqn_appr_rate_main}. Allowing the coefficient norm to grow like $n^\alpha$ weakens the corresponding regularity and approximation order by the same factor.

Our first result identifies the exact approximation space at this scale. For $0<\alpha<\frac{d+2k+1}{2d}$ and every quasi-uniform sequence $\{\Theta_n\}_{n=1}^\infty$,
\begin{equation}\label{eqn_shallowk_bound}
    f\in\mathcal B_{2,\infty,*}^{\frac{d+2k+1}{2}-\alpha d}(\SS^d)\quad\Longleftrightarrow\quad\text{there exists }M>0\text{ such that }\inf_{f_n\in L_{n,M}^{k,\alpha}(\Theta_n)}\|f-f_n\|_{\mathcal L^2(\SS^d)}\lesssim n^{-\frac{d+2k+1}{2d}+\alpha}.
\end{equation}
Thus there is no loss of smoothness in the inverse conclusion. The rate assumption remains essential: when $\alpha>0$, mere convergence by approximants satisfying the same coefficient growth bound does not imply membership in $\mathcal B_{2,\infty,*}^{\frac{d+2k+1}{2}-\alpha d}(\SS^d)$. This contrasts with the endpoint case $\alpha=0$ in \eqref{eqn:appr_equiv_sob}.

We next remove the coefficient constraint and retain only the approximation rate. The resulting picture is more delicate. If $0<\alpha<\frac{1}{2}$ and $\{\Theta_n\}_{n=1}^\infty$ is nested, the estimate
\begin{equation}\label{eqn_rate_norestr}
    \inf_{f_n\in L_n^k(\Theta_n)}\|f-f_n\|_{\mathcal L^2(\SS^d)}\lesssim n^{-\frac{d+2k+1}{2d}+\alpha}
\end{equation}
does not imply that $f\in\mathcal B_{2,\infty,*}^{\frac{d+2k+1}{2}-\alpha d}(\SS^d)$. Indeed, nestedness preserves at least one ridge function at all levels, and a single ridge has exact Besov smoothness $k+\frac{1}{2}$. If $\alpha>\frac{1}{2}$ and the nested sequence is antipodally quasi-uniform, the Bernstein inequality yields the converse implication at the endpoint order $\frac{d+2k+1}{2}-\alpha d$.

The inverse conclusion can nevertheless be recovered from two geometrically independent parameter sequences. Suppose the same target satisfies \eqref{eqn_rate_norestr} for quasi-uniform sequences $\{\Theta_n\}_{n=1}^\infty$ and $\{\widetilde\Theta_n\}_{n=1}^\infty$, and that $\Theta_n\cup\widetilde\Theta_n$ is antipodally separated at scale $n^{-\frac{1}{d}}$. Subtracting approximants on the two parameter sets and applying a coefficient lower bound on their union yields
\[
    \sul_{j=1}^na_{j,n}^2\lesssim n^{2\alpha-1}.
\]
Thus the coefficient growth required above is recovered from the two approximation estimates themselves, and the Besov inverse conclusion follows from \eqref{eqn_shallowk_bound}.

These results separate two pieces of information that are often conflated in neural-network approximation: the size of the approximation error and the stability of the representation producing it. A single parameter sequence may approximate a target rapidly for geometric reasons that reveal little about its smoothness. An explicit coefficient bound, or approximation on two suitably separated parameter sequences, rules out this degeneracy and restores a sharp inverse statement. The result therefore has a saturation-type interpretation once the admissible coefficient growth is fixed, but it is not a classical saturation theorem for a single fixed approximation process. This perspective complements recent work relating neural-network norms, approximation classes, Bernstein inequalities, and generalization \cite{siegel2022sharp,siegel2022optimal,siegel2023characterization,mhaskar2023tractability,huang2023generalization}.

The remainder of the paper is organized as follows. Section~2 states the main theorem. Section~3 introduces the spherical harmonic notation, the ReLU$^k$ expansion, and the quadrature formula used in the analysis. Section~4 contains the proofs.

\section{Main Results}

\begin{theorem}\label{thm:main}
Let $d\in\NN$, $k\in\NN_0$, $0<\alpha<\frac{d+2k+1}{2d}$, and let $\{\Theta_n\}_{n=1}^\infty$ be a sequence of quasi-uniform collections on $\SS^d$.
\begin{enumerate}
\item\label{item1} Let $f\in\mathcal L^2(\SS^d)$, then $f\in\mathcal B_{2,\infty,*}^{\frac{d+2k+1}{2}-\alpha d}(\SS^d)$ if and only if there exists $M$ independent of $n$, such that
\begin{equation}\label{eqn:item1_rate}
\inf_{f_n\in L_{n,M}^{k,\alpha}(\Theta_n)}\|f-f_n\|_{\mathcal L^2(\SS^d)}\lesssim n^{-\frac{d+2k+1}{2d}+\alpha}.
\end{equation}
\item\label{item2} For every $M>0$ and $f\in\mathcal L_*^2(\SS^d)$,
\begin{equation}\label{eqn:coeff_alpha}
\lim_{n\to\infty}
\inf_{f_n\in L_{n,M}^{k,\alpha}(\Theta_n)}
\|f-f_n\|_{\mathcal L^2(\SS^d)}=0.
\end{equation}
\item\label{item3} Suppose further that $\{\Theta_n\}_{n=1}^\infty$ is antipodally quasi-uniform and nested. If $\alpha>\frac12$ and $f\in\mathcal L^2(\SS^d)$ satisfies
\begin{equation}\label{eqn:rate_alpha}
\inf_{f_n\in L_n^k(\Theta_n)}\|f-f_n\|_{\mathcal L^2(\SS^d)}\lesssim n^{-\frac{d+2k+1}{2d}+\alpha},
\end{equation}
then it follows that $f\in \mathcal B_{2,\infty,*}^{\frac{d+2k+1}{2}-\alpha d}(\SS^d)$.
On the other hand, for every $\alpha \leq 1/2$, this characterization fails, i.e., there exists an $f\notin \mathcal B_{2,\infty,*}^{\frac{d+2k+1}{2}-\alpha d}(\SS^d)$ which satisfies \eqref{eqn:rate_alpha}.
\item\label{item4} Let $f\in\mathcal L^2(\SS^d)$, and suppose there exists another sequence of antipodally quasi-uniform collections $\{\widetilde\Theta_n\}_{n=1}^\infty$, with $\widetilde\Theta_n=\{\widetilde\theta_{j,n}^*\}_{j=1}^n$, such that $\#(\Theta_n\cup\widetilde\Theta_n)=2n$,
\begin{equation}\label{eqn:item4_union_separation}
\min_{\substack{\xi,\zeta\in\Theta_n\cup\widetilde\Theta_n\\\xi\neq\zeta}}\min\{\rho(\xi,\zeta),\rho(\xi,-\zeta)\}\gtrsim n^{-\frac{1}{d}},
\end{equation}
and
\begin{equation}\label{eqn:item4_two_rates}
\inf_{f_n\in L_n^k(\Theta_n)}\|f-f_n\|_{\mathcal L^2(\SS^d)}+\inf_{\widetilde f_n\in L_n^k(\widetilde\Theta_n)}\|f-\widetilde f_n\|_{\mathcal L^2(\SS^d)}\lesssim n^{-\frac{d+2k+1}{2d}+\alpha}.
\end{equation}
Then $f\in\mathcal B_{2,\infty,*}^{\frac{d+2k+1}{2}-\alpha d}(\SS^d)$.
\end{enumerate}
\end{theorem}

\section{Preliminaries}

Throughout, $d\eta$ denotes the normalized rotation-invariant measure on $\SS^d$. We write $\PP_m(\SS^d)$ for the restrictions to $\SS^d$ of polynomials on $\RR^{d+1}$ of total degree at most $m$.

\subsection{Spherical harmonics, Besov spaces, and Sobolev spaces}

For $m\in\NN_0$, let $\mathcal Y_m$ be the space of spherical harmonics of degree $m$, and fix an orthonormal basis $\{Y_{m,\ell}\}_{\ell=1}^{N(m)}$. Its dimension is
\begin{equation}
    N(m):=\dim(\mathcal Y_m)=\dim(\PP_m(\SS^d))-\dim(\PP_{m-1}(\SS^d))=\begin{cases}
        \displaystyle(2m+d-1)\frac{\Gamma(m+d-1)}{\Gamma(m+1)\Gamma(d)},&m\ne0,\\
        1,&m=0.
    \end{cases}
\end{equation}
The orthogonal decomposition into spherical harmonics gives, for every $f\in\mathcal L^2(\SS^d)$,
\begin{equation}
    f(\eta)=\sum_{m=0}^\infty\sum_{\ell=1}^{N(m)}\widehat f(m,\ell)Y_{m,\ell}(\eta),\qquad\text{a.e. }\eta\in\SS^d,
\end{equation}
where
\[
    \widehat f(m,\ell)=\langle f,Y_{m,\ell}\rangle_{\mathcal L^2(\SS^d)},\qquad m\in\NN_0,\quad1\leq\ell\leq N(m).
\]

Let $\Pi_m$ be the $\mathcal L^2(\SS^d)$-orthogonal projection onto $\mathcal Y_m$, and for $J\in\NN_0$ set
\begin{equation}\label{eqn:def_sharp_projectors}
    \mathcal P_Jf:=\sum_{m=0}^J\Pi_mf,\qquad \mathcal P_J^cf:=f-\mathcal P_Jf,
\end{equation}
as well as
\begin{equation}\label{eqn:def_dyadic_blocks}
    \Delta_0f:=\mathcal P_1f,\qquad \Delta_qf:=(\mathcal P_{2^q}-\mathcal P_{2^{q-1}})f,\quad q\geq1.
\end{equation}

\begin{definition}[Besov and Sobolev spaces]\label{def_besov_sphere}
For $\gamma>0$, define
\[
    \mathcal B_{2,\infty}^\gamma(\SS^d):=\Big\{f\in\mathcal L^2(\SS^d):~\|f\|_{\mathcal B_{2,\infty}^\gamma(\SS^d)}<\infty\Big\},
\]
where
\begin{equation}\label{eqn:def_besov_norm}
    \|f\|_{\mathcal B_{2,\infty}^\gamma(\SS^d)}:=\|\Delta_0f\|_{\mathcal L^2(\SS^d)}+\sup_{q\geq1}2^{q\gamma}\|\Delta_qf\|_{\mathcal L^2(\SS^d)}.
\end{equation}
For $\gamma\geq0$, the Sobolev space $\mathcal H^\gamma(\SS^d)$ is defined by
\begin{equation}\label{eqn:Sob_norm_Parseval}
    \|f\|_{\mathcal H^\gamma(\SS^d)}^2:=\sum_{m=0}^\infty(1+m)^{2\gamma}\sum_{\ell=1}^{N(m)}|\widehat f(m,\ell)|^2.
\end{equation}
For $\gamma=0$, this is the norm on $\mathcal L^2(\SS^d)$.
\end{definition}

\subsection{Quasi-uniform parameter collections}

Let $\rho$ be the geodesic distance on $\SS^d$. For $\Theta_n=\{\theta_{j,n}^*\}_{j=1}^n\subset\SS^d$, define its fill distance, separation distance, and antipodal separation distance by
\begin{equation*}
    h_n:=\mal_{\eta\in\SS^d}\mil_{1\leq j\leq n}\rho(\eta,\theta_{j,n}^*),\qquad\widetilde h_n:=\mil_{i\ne j}\rho(\theta_{i,n}^*,\theta_{j,n}^*),
\end{equation*}
and
\begin{equation*}
    \underline h_n:=\min\Big\{\mil_{i\ne j}\rho(\theta_{i,n}^*,\theta_{j,n}^*),\mil_{i\ne j}\rho(-\theta_{i,n}^*,\theta_{j,n}^*)\Big\}.
\end{equation*}
The collection $\Theta_n$ is called quasi-uniform if $h_n\lesssim\widetilde h_n$, and antipodally quasi-uniform if $h_n\lesssim\underline h_n$. A sequence $\{\Theta_n\}_{n=1}^\infty$ has either property when the corresponding constant is uniform in $n$, and it is called nested if $\Theta_n\subset\Theta_{n+1}$ for every $n$. In particular,
\[
    h_n\simeq\widetilde h_n\simeq n^{-\frac{1}{d}}
\]
for a quasi-uniform sequence, while $h_n\simeq\underline h_n\simeq n^{-\frac{1}{d}}$ in the antipodally quasi-uniform case. When one value of $n$ is fixed, we write $h=h_n$.

\subsection{Ultraspherical expansion of \texorpdfstring{$\sigma_k$}{sigma-k}}

Let $\mathcal L^2_{w_d}([-1,1])$ be the weighted space with
\[
    \langle f,g\rangle_{w_d}=\int_{-1}^1f(t)g(t)(1-t^2)^{\frac{d-2}{2}}dt,\qquad\|f\|_{\mathcal L^2_{w_d}([-1,1])}=\langle f,f\rangle_{w_d}^{\frac{1}{2}}.
\]
The ultraspherical polynomials $\{p_m\}_{m=0}^\infty$ form an orthogonal basis of this space and satisfy
\begin{equation}\label{eqn:def_ultraspherical}
    \langle p_m,p_\ell\rangle_{w_d}=0,\qquad m\ne\ell,\quad m,\ell\in\NN_0.
\end{equation}
We use the normalization for which the spherical addition formula reads
\begin{equation}\label{eqn:sum_Y_nl}
    p_m(\eta\cdot\theta)=\sum_{\ell=1}^{N(m)}Y_{m,\ell}(\eta)Y_{m,\ell}(\theta);
\end{equation}
see \cite[Theorem~1.2.6]{dai2013approximation}. Here $\omega_d$ denotes the surface area of $\SS^d$.

The activation $\sigma_k\in\mathcal L^2_{w_d}([-1,1])$ admits the expansion
\begin{equation}\label{eqn:ultraspherical_expansion_ReLUk}
    \sigma_k=\sum_{m=0}^\infty\widehat{\sigma_k}(m)p_m,
\end{equation}
where
\[
    \widehat{\sigma_k}(m)=\frac{\langle p_m,\sigma_k\rangle_{w_d}}{\|p_m\|_{\mathcal L^2_{w_d}([-1,1])}^2}.
\]
Following \cite{schneider1967problem,bourgain2006projection,mhaskar2006weighted,bach2017breaking}, the nonzero coefficients are indexed by
\[
    E_{\sigma_k}=\{m\geq k+1:m-k\text{ odd}\}\cup\{0,\ldots,k\}.
\]
For $m\geq k+1$ with $m-k$ odd, they are given explicitly by
\begin{equation}\label{eqn:hat_sigma_large}
    \widehat{\sigma_k}(m)=\frac{\omega_{d-1}k!\Gamma(\frac{d}{2})}{\omega_d}\frac{(-1)^{\frac{m-k-1}{2}}\Gamma(m-k)}{2^m\Gamma\!\Big(\frac{m-k+1}{2}\Big)\Gamma\!\Big(\frac{m+d+k+1}{2}\Big)}.
\end{equation}
Then a neuron on $\SS^d$ has an expansion
\begin{equation}\label{eqn_neu_expan}
    \sigma_k(\theta\cdot\eta)=\sul_{m\in E_{\sigma_k}}\widehat{\sigma_k}(m)\sul_{\ell=1}^{N(m)}Y_{m,\ell}(\theta)Y_{m,\ell}(\eta).
\end{equation}

Related ultraspherical expansions of ReLU-type zonal functions appear in \cite{mhaskar2019reluZonal}.

\subsection{Polynomial quadrature on quasi-uniform collections}

We shall use the following positive quadrature rule.

\begin{lemma}[{\cite[Theorem~4.1]{mhaskar2001spherical}, \cite[Theorem~3.2]{narcowich2006sobolev}}]\label{lem:quadrature}
There exist nonnegative weights $\tau_{1,n},\dots,\tau_{n,n}$ satisfying $\tau_{j,n}\lesssim h^d$, and a constant $c_1$ independent of $n$ and $h$, such that
\begin{equation}\label{eqn:quadrature_formula}
    \int_{\SS^d}p(\eta)q(\eta)d\eta=\sul_{j=1}^n\tau_{j,n}p(\theta_{j,n}^*)q(\theta_{j,n}^*),\qquad\forall p,q\in\PP_J(\SS^d),
\end{equation}
where $J=\lfloor c_1h^{-1}\rfloor$. Moreover, if
\begin{equation}
    \mil_{i\ne j}\rho(\theta_{i,n}^*,\theta_{j,n}^*)\gtrsim h,
\end{equation}
then there exists $c_2>0$ such that \eqref{eqn:quadrature_formula} holds for $J\leq c_2n^{\frac{1}{d}}$ and
\begin{equation}
    \tau_{j,n}\simeq h^d\simeq n^{-1},\qquad j=1,\ldots,n.
\end{equation}
The constants are independent of $n$ and may depend only on the dimension and the uniform quasi-uniformity constants.
\end{lemma}

We shall use the endpoint case of the direct theorem in \cite{liu2025integral} in the following form.

\begin{theorem}[Theorem 4.1, \cite{liu2025integral}]\label{thm:endpoint_direct}
Let $d\geq2$, $k\in\NN_0$, and $\{\Theta_n\}_{n=1}^\infty$ be quasi-uniform. For every $f\in\mathcal H_*^{\frac{d+2k+1}{2}}(\SS^d)$, there is
\[
    f_n=\sum_{j=1}^na_{j,n}\sigma_k(\theta_{j,n}^*\cdot\circ)
\]
such that
\begin{equation}\label{eqn:endpoint_direct}
    \|f-f_n\|_{\mathcal L^2(\SS^d)}\lesssim n^{-\frac{d+2k+1}{2d}}\|f\|_{\mathcal H^{\frac{d+2k+1}{2}}(\SS^d)},\qquad \Big(n\sum_{j=1}^na_{j,n}^2\Big)^{\frac{1}{2}}\lesssim\|f\|_{\mathcal H^{\frac{d+2k+1}{2}}(\SS^d)}.
\end{equation}
\end{theorem}

The theorem is the endpoint specialization of the spherical direct approximation theorem in \cite{liu2025integral}. Although its parity condition is commonly stated for the pure even or odd component, the same construction covers the additional $\PP_k(\SS^d)$ component: all degrees $0\leq m\leq k$ belong to $E_{\sigma_k}$ and are reproduced by the quadrature formula.

We also record the stability estimate used in items~\ref{item3} and \ref{item4}. Its Bernstein part is the estimate proved in \cite{liu2026optimal}; the same localization argument also gives the coefficient lower bound. We include the short common argument so that the latter is not left implicit.

\begin{theorem}[Theorem 1, \cite{liu2026optimal}; Theorem 1, \cite{liu2025condition}]\label{lem:ridge_stability}
Let $d\geq2$, $k\in\NN_0$, and $\{\Theta_n\}_{n=1}^\infty$ be antipodally quasi-uniform. Then
\begin{equation}\label{eqn:coefficient_lower_general}
    \inf_{a_n\in\RR^n,~\|a_n\|_2=1}\Big\|\sum_{j=1}^na_{j,n}\sigma_k(\theta_{j,n}^*\cdot\circ)\Big\|_{\mathcal L^2(\SS^d)}\simeq n^{-\frac{2k+1}{2d}}.
\end{equation}
Moreover, let $0\leq s < r<k+\frac{1}{2}$, then for $f_n\in L_n^k(\Theta_n)$,
\begin{equation}\label{eqn:bernstein_general}
    \|f_n\|_{\mathcal{H}^r(\SS^d)}\lesssim n^{\frac{r-s}{d}}\|f_n\|_{\mathcal{H}^s(\SS^d)},
\end{equation}
where the implied constant is independent of $n$, $\Theta_n$, and $f_n$.
\end{theorem}

\begin{proof}
    We write
    \begin{equation*}
        \begin{split}
            &\Big\|\sum_{j=1}^na_{j,n}\sigma_k(\theta_{j,n}^*\cdot\circ)\Big\|_{\mathcal L^2(\SS^d)}^2=\int_{\SS^d}\Big(\sum_{j=1}^na_{j,n}\sigma_k(\theta_{j,n}^*\cdot\eta)\Big)^2d\eta\\
            =&\int_{\SS^d}a_n^\top\Big(\sigma_k(\theta_{i,n}^*\cdot\eta)\sigma_k(\theta_{j,n}^*\cdot\eta)\Big)_{i,j=1}^na_nd\eta=a_n^\top\mathcal M_na_n,
        \end{split}
    \end{equation*}
    where
    \begin{equation*}
        \mathcal M_n=\Bigg(\int_{\SS^d}\sigma_k(\theta_{i,n}^*\cdot\eta)\sigma_k(\theta_{j,n}^*\cdot\eta)d\eta\Bigg)_{i,j=1}^n.
    \end{equation*}
    By \cite[Theorem 1]{liu2025condition},
    $$\lambda_{\min}(\mathcal M_n)\simeq n^{-\frac{2k+1}{d}},$$
    which proves \eqref{eqn:coefficient_lower_general}.

    \eqref{eqn:bernstein_general} is exactly the result in \cite{liu2026optimal}
\end{proof}

\section{Proofs of the main theorem}
Throughout this section, $c_2$ and $\tau_{j,n}$ are as in Lemma~\ref{lem:quadrature}. In particular, $\tau_{j,n}\simeq n^{-1}$.

\subsection{Proof of item \ref{item1}}
\paragraph{Direct implication.}
Let $f\in\mathcal B_{2,\infty,*}^{\frac{d+2k+1}{2}-\alpha d}(\SS^d)$. Choose $Q=Q(n)$ so that $2^Q\leq n^{\frac{1}{d}}<2^{Q+1}$ and put $g_n:=\mathcal P_{2^Q}f$. The dyadic definition of the Besov norm and the orthogonality of the blocks give
\begin{equation}\label{eqn:item1_truncation_bounds}
    \|f-g_n\|_{\mathcal L^2(\SS^d)}\lesssim2^{-Q(\frac{d+2k+1}{2}-\alpha d)}\|f\|_{\mathcal B_{2,\infty}^{\frac{d+2k+1}{2}-\alpha d}(\SS^d)},\qquad \|g_n\|_{\mathcal H^{\frac{d+2k+1}{2}}(\SS^d)}\lesssim2^{Q\alpha d}\|f\|_{\mathcal B_{2,\infty}^{\frac{d+2k+1}{2}-\alpha d}(\SS^d)}.
\end{equation}
Indeed, the first estimate follows by summing the squared high-frequency blocks, while the second follows from
\[
    \|g_n\|_{\mathcal H^{\frac{d+2k+1}{2}}}^2\lesssim\sum_{q=0}^Q2^{q(d+2k+1)}\|\Delta_qf\|_2^2\lesssim\|f\|_{\mathcal B_{2,\infty}^{\frac{d+2k+1}{2}-\alpha d}}^2\sum_{q=0}^Q2^{2q(\alpha d)}\lesssim2^{2Q\alpha d}\|f\|_{\mathcal B_{2,\infty}^{\frac{d+2k+1}{2}-\alpha d}}^2,
\]
where $\alpha d>0$. Apply Theorem~\ref{thm:endpoint_direct} to $g_n$. The resulting network, denoted by $f_n$, satisfies
\[
    \Big(n\sum_{j=1}^na_{j,n}^2\Big)^{\frac{1}{2}}\lesssim\|g_n\|_{\mathcal H^{\frac{d+2k+1}{2}}}\lesssim n^\alpha\|f\|_{\mathcal B_{2,\infty}^{\frac{d+2k+1}{2}-\alpha d}},
\]
and hence belongs to $L_{n,M}^{k,\alpha}(\Theta_n)$ for an $M$ independent of $n$. Moreover,
\begin{equation*}
    \begin{split}
        \|f-f_n\|_2\leq&\|f-g_n\|_2+\|g_n-f_n\|_2\lesssim2^{-Q(\frac{d+2k+1}{2}-\alpha d)}\|f\|_{\mathcal B_{2,\infty}^{\frac{d+2k+1}{2}-\alpha d}}+n^{-\frac{d+2k+1}{2d}}2^{Q\alpha d}\|f\|_{\mathcal B_{2,\infty}^{\frac{d+2k+1}{2}-\alpha d}}\\
        \lesssim& n^{-\frac{d+2k+1}{2d}+\alpha}\|f\|_{\mathcal B_{2,\infty}^{\frac{d+2k+1}{2}-\alpha d}}.
    \end{split}
\end{equation*}

\paragraph{Inverse implication.}
For each $n$, let $J=\lfloor c_2n^{\frac{1}{d}}\rfloor$.
For $0\le m\le J$ and $1\le\ell\le N(m)$, define
\begin{equation}\label{eqn_def_Vml}
    V_n(m,\ell):=\bigl(\sqrt{\tau_{1,n}}Y_{m,\ell}(\theta_{1,n}^*),\ldots,\sqrt{\tau_{n,n}}Y_{m,\ell}(\theta_{n,n}^*)\bigr)^\top.
\end{equation}
Then
\begin{equation}\label{eqn:V_orthonormal}
    V_n(m,\ell)^\top V_n(m',\ell')=\delta_{(m,\ell),(m',\ell')},\qquad 0\le m,m'\le J,
\end{equation}
by applying \eqref{eqn:quadrature_formula} to $Y_{m,\ell}Y_{m',\ell'}$.

For a representation $f_n=\sum_{j=1}^na_{j,n}\sigma_k(\theta_{j,n}^*\cdot\circ)$ satisfying the defining constraint of $L_{n,M}^{k,\alpha}$, set
\begin{equation}\label{eqn:y_v_def}
    y_n:=\Big(\frac{a_{1,n}}{\sqrt{\tau_{1,n}}},\ldots,\frac{a_{n,n}}{\sqrt{\tau_{n,n}}}\Big)^\top,\qquad \widehat{v_n}(m,\ell):=y_n^\top V_n(m,\ell).
\end{equation}
For $m\le J$, the addition formula gives
\begin{equation}\label{eqn:fn_harmonic_coeff}
    \widehat f_n(m,\ell)=\widehat{\sigma_k}(m)\widehat{v_n}(m,\ell).
\end{equation}
Moreover, by \eqref{eqn:V_orthonormal}, Bessel's inequality, and $\tau_{j,n}\simeq n^{-1}$,
\begin{equation}\label{eqn:y_bound}
    \sum_{m=0}^{J}\sum_{\ell=1}^{N(m)}|\widehat{v_n}(m,\ell)|^2
    \le\|y_n\|_2^2
    =\sum_{j=1}^n\frac{a_{j,n}^2}{\tau_{j,n}}
    \lesssim n\sum_{j=1}^na_{j,n}^2.
\end{equation}

For every sufficiently large $q$, choose $n=n_q:=\lceil C2^{qd}\rceil$, where the fixed constant $C$ is large enough that $2^q\leq\lfloor c_2n_q^{\frac{1}{d}}\rfloor=J$. Thus $n_q\simeq2^{qd}$. By the assumed approximation estimate, there exists $f_n\in L_{n,M}^{k,\alpha}$ such that
\begin{equation}\label{eqn:item1_inverse_approx_nq}
    \|f-f_n\|_{\mathcal L^2(\SS^d)}\lesssim n^{-\frac{d+2k+1}{2d}+\alpha}\lesssim2^{-q(\frac{d+2k+1}{2}-\alpha d)}.
\end{equation}
The explicit formula \eqref{eqn:hat_sigma_large} and the Gamma-ratio asymptotics imply
\begin{equation}\label{eqn:sigma_decay_r}
    |\widehat{\sigma_k}(m)|\lesssim m^{-\frac{d+2k+1}{2}},\qquad m>k,
\end{equation}
where $\widehat{\sigma_k}(m)=0$ for $m>k$ with $m-k$ even. Hence, for all sufficiently large $q$, \eqref{eqn:fn_harmonic_coeff} and \eqref{eqn:y_bound} give
\begin{equation}\label{eqn:item1_inverse_shell_fn}
    \begin{split}
        \|(\mathcal P_{2^q}-\mathcal P_{2^{q-1}})f_n\|_{\mathcal L^2(\SS^d)}^2
        &=\sum_{m=2^{q-1}+1}^{2^q}\widehat{\sigma_k}(m)^2
        \sum_{\ell=1}^{N(m)}|\widehat{v_n}(m,\ell)|^2\\
        &\le\max_{2^{q-1}<m\le2^q}|\widehat{\sigma_k}(m)|^2
        \|y_n\|_2^2\\
        &\lesssim2^{-q(d+2k+1)}n^{2\alpha}\lesssim2^{-2q(\frac{d+2k+1}{2}-\alpha d)}.
    \end{split}
\end{equation}
Since $\mathcal P_{2^q}-\mathcal P_{2^{q-1}}$ is an orthogonal projection, \eqref{eqn:item1_inverse_approx_nq} and \eqref{eqn:item1_inverse_shell_fn} yield
\begin{equation}\label{eqn:item1_inverse_shell_f}
    \begin{split}
        \|(\mathcal P_{2^q}-\mathcal P_{2^{q-1}})f\|_{\mathcal L^2(\SS^d)}^2
        &\le2\|(\mathcal P_{2^q}-\mathcal P_{2^{q-1}})f_n\|_{\mathcal L^2(\SS^d)}^2+2\|(\mathcal P_{2^q}-\mathcal P_{2^{q-1}})(f-f_n)\|_{\mathcal L^2(\SS^d)}^2\\
        &\lesssim2^{-2q(\frac{d+2k+1}{2}-\alpha d)}.
    \end{split}
\end{equation}

Equation \eqref{eqn:item1_inverse_shell_f} and the finitely many low-frequency blocks give
\[
    \sup_{q\geq0}2^{q(\frac{d+2k+1}{2}-\alpha d)}\|\Delta_qf\|_{\mathcal L^2(\SS^d)}<\infty,
\]
so $f\in\mathcal B_{2,\infty}^{\frac{d+2k+1}{2}-\alpha d}(\SS^d)$. It remains to verify the parity condition. For every $m>k$ with $m-k$ even, one has $\widehat{\sigma_k}(m)=0$, so $\Pi_mf_n=0$. Since $\frac{d+2k+1}{2}-\alpha d>0$, \eqref{eqn:item1_inverse_approx_nq} implies $f_n\to f$ in $\mathcal L^2(\SS^d)$ along the above choice of $n$. Consequently, $\Pi_mf=0$ for every such $m$, and hence $f\in\mathcal B_{2,\infty,*}^{\frac{d+2k+1}{2}-\alpha d}(\SS^d)$.

\subsection{Proof of item \ref{item2}}
    Fix arbitrary $\epsilon>0$ and fix $L$ such that
    \begin{equation}
        \|f-\Pi_Lf\|_{\mathcal{L}^2(\SS^d)}\leq\epsilon.
    \end{equation}
    By Theorem \ref{thm:endpoint_direct}, for each $n\in\NN$, there exists $M_0>0$ and $a_{n,1},\dots,a_{n,n}$ with
    \begin{equation*}
        n\sul_{j=1}^na_{j,n}^2\leq M_0\|\Pi_Lf\|_{\mathcal{H}^{\frac{d+2k+1}{2}}(\SS^d)}
    \end{equation*}
    such that
    \begin{equation*}
        \|f_n-\Pi_Lf\|_{\mathcal{L}^2(\SS^d)}\lesssim\|\Pi_Lf\|_{\mathcal{H}^{\frac{d+2k+1}{2}}(\SS^d)}n^{-\frac{d+2k+1}{2d}}.
    \end{equation*}
    Let $n$ be sufficiently large such that
    $$Mn^\alpha>M_0\|\Pi_Lf\|_{\mathcal{H}^{\frac{d+2k+1}{2}}(\SS^d)},\qquad\|\Pi_Lf\|_{\mathcal{H}^{\frac{d+2k+1}{2}}(\SS^d)}n^{-\frac{d+2k+1}{2d}}\leq\epsilon,$$
    then $f_n\in L_{n,M}^{k,\alpha}(\Theta_n)$ and
    \begin{equation}
        \|f-f_n\|_{\mathcal{L}^2(\SS^d)}\leq\|f-\Pi_Lf\|_{\mathcal{L}^2(\SS^d)}+\|f_n-\Pi_Lf\|_{\mathcal{L}^2(\SS^d)}\leq2\epsilon.
    \end{equation}

\subsection{Proof of item \ref{item3}}
\paragraph{Case 1: $\alpha<\frac{1}{2}$.}
Since $\{\Theta_n\}_{n=1}^\infty$ is nested, we may fix $\theta_{1,1}^*\in\Theta_1$, so that $\theta_{1,1}^*\in\Theta_n$ for every $n$. Let
\[
    f(\eta):=\sigma_k(\theta_{1,1}^*\cdot\eta).
\]
Then $f\in L_n^k(\Theta_n)$ for every $n$, so the left-hand side of \eqref{eqn:rate_alpha} is identically zero.

It remains to verify the claimed lack of regularity. By \eqref{eqn:sum_Y_nl}, \eqref{eqn:ultraspherical_expansion_ReLUk}, and \eqref{eqn:hat_sigma_large}, for $m\ge k+1$ with $m-k$ odd,
\begin{equation}\label{eqn:item3_harmonic_size}
    \|\Pi_mf\|_{\mathcal L^2(\SS^d)}^2=\widehat{\sigma_k}(m)^2\sum_{\ell=1}^{N(m)}|Y_{m,\ell}(\theta_{1,1}^*)|^2=\widehat{\sigma_k}(m)^2N(m)\asymp m^{-2k-2}.
\end{equation}
Consequently, for every sufficiently large $q$,
\[
    \|\Delta_qf\|_{\mathcal L^2(\SS^d)}^2\asymp\sum_{\substack{2^{q-1}<m\leq2^q\\m-k\ \mathrm{odd}}}m^{-2k-2}\asymp2^{-q(2k+1)}.
\]
Therefore,
\[
    2^{q(\frac{d+2k+1}{2}-\alpha d)}\|\Delta_qf\|_{\mathcal L^2(\SS^d)}\asymp2^{q(\frac{d+2k+1}{2}-\alpha d-k-\frac{1}{2})}=2^{qd(\frac{1}{2}-\alpha)}\longrightarrow\infty.
\]
The expansion of $f$ has the required parity above degree $k$, and hence $f\in\mathcal L_*^2(\SS^d)\setminus\mathcal B_{2,\infty,*}^{\frac{d+2k+1}{2}-\alpha d}(\SS^d)$.

\paragraph{Case 2: $\alpha>\frac{1}{2}$.}
Let $f\in\mathcal L^2(\SS^d)$ satisfy \eqref{eqn:rate_alpha}, for $q\geq0$, let $n_q=2^{(q+1)d}$ and choose $v_q\in L_{n_q}^k(\Theta_{n_q})$ such that
\begin{equation}\label{eqn:item3_dyadic_approximants}
    \|f-v_q\|_{\mathcal L^2(\SS^d)}\lesssim2^{-q(\frac{d+2k+1}{2}-\alpha d)}.
\end{equation}
Set $u_0:=v_0$ and $u_q:=v_q-v_{q-1}$ for $q\geq1$. By nestedness, $u_q\in L_{n_q}^k(\Theta_{n_q})$, and \eqref{eqn:item3_dyadic_approximants} gives $\|u_q\|_{\mathcal L^2(\SS^d)}\lesssim2^{-q(\frac{d+2k+1}{2}-\alpha d)}$ after adjusting the first term. Take $r\in(\frac{d+2k+1}{2}-\alpha d,k+\frac{1}{2})$, then the antipodal quasi-uniformity assumption and \eqref{eqn:bernstein_general} yield
\begin{equation}\label{eqn:item3_difference_bernstein}
    \|u_q\|_{\mathcal H^r(\SS^d)}\lesssim n_q^{\frac{r}{d}}\|u_q\|_{\mathcal L^2(\SS^d)}\lesssim2^{q(r-\frac{d+2k+1}{2}+\alpha d)}.
\end{equation}
As $f=\sul_{q=0}^\infty u_q$ in $\mathcal L^2(\SS^d)$, for $Q\in\NN$, equations \eqref{eqn:item3_dyadic_approximants}--\eqref{eqn:item3_difference_bernstein} give
\begin{equation}\label{eqn:item3_besov_block}
    \begin{split}
        \|\Delta_Qf\|_{\mathcal L^2(\SS^d)}&\leq\sum_{q<Q}\|\Delta_Qu_q\|_{\mathcal L^2(\SS^d)}+\sum_{q\geq Q}\|u_q\|_{\mathcal L^2(\SS^d)}\lesssim2^{-Qr}\sum_{q<Q}\|\Delta_Qu_q\|_{\mathcal H^r(\SS^d)}+\sum_{q\geq Q}\|u_q\|_{\mathcal L^2(\SS^d)}\\
        &\lesssim2^{-Qr}\sum_{q<Q}2^{q(r-\frac{d+2k+1}{2}+\alpha d)}+\sum_{q\geq Q}2^{-q(\frac{d+2k+1}{2}-\alpha d)}\lesssim2^{-Q(\frac{d+2k+1}{2}-\alpha d)}.
    \end{split}
\end{equation}
Hence $f\in\mathcal B_{2,\infty}^{\frac{d+2k+1}{2}-\alpha d}(\SS^d)$. Finally, $v_q\to f$ in $\mathcal L^2$, and every $v_q$ has zero projection onto $\mathcal Y_m$ whenever $m>k$ and $m-k$ is even. Therefore $f\in\mathcal B_{2,\infty,*}^{\frac{d+2k+1}{2}-\alpha d}(\SS^d)$.

\paragraph{Case 3: $\alpha=\frac{1}{2}$.}
By nestedness, after relabeling the points if necessary, we may write
\[
\Theta_n=\{\theta_1^*,\ldots,\theta_n^*\},\qquad n=1,2\dots.
\]
Since $\dim\PP_J(\SS^d)\asymp J^d$, we may choose $s_0$ sufficiently large and integers $J_s$, $s\geq s_0$, such that
\begin{equation}\label{eqn:item3_endpoint_Js}
J_s+2\in E_{\sigma_k},\qquad J_s\asymp2^{\frac{s+1}{d}},\qquad \dim\PP_{J_s}(\SS^d)\leq2^{s-1},
\end{equation}
and $J_s+2\leq c_22^{\frac{s+1}{d}}$, which lies below the quadrature degree associated with $\Theta_{2^{s+1}}$.

For each $s\geq s_0$, we follow the notations in \cite[(4.28)--(4.34)]{liu2025condition} and denote the orthogonal complement of the orthonormal vectors in \eqref{eqn:V_orthonormal}:
\[
\mathcal V(J_s):=\operatorname{span}\left\{V_{2^{s+1}}(m,\ell):m\in E_{\sigma_k},~m\leq J_s,\ 1\leq\ell\leq N(m)\right\}^{\perp}.
\]

By \eqref{eqn:V_orthonormal} and \eqref{eqn:item3_endpoint_Js},
\[
\begin{split}
\dim\left(\mathcal V(J_s)\cap\left\{v\in\RR^{2^{s+1}}:v_1=\cdots=v_{2^s}=0\right\}\right)
\geq2^s-\dim\PP_{J_s}(\SS^d)\geq2^{s-1}.
\end{split}
\]
We may therefore choose a unit vector
\[
v^{(s)}=\bigl(v_1^{(s)},\ldots,v_{2^{s+1}}^{(s)}\bigr)^\top\in\mathcal V(J_s),\qquad v_1^{(s)}=\cdots=v_{2^s}^{(s)}=0.
\]
Thus the coefficients of the corresponding ridge combination are supported precisely on
\[
\Theta_{2^{s+1}}\setminus\Theta_{2^s}=\{\theta_{2^s+1}^*,\ldots,\theta_{2^{s+1}}^*\}.
\]

Let
\[
W_s:=\operatorname{diag}\bigl(\tau_{1,2^{s+1}},\ldots,\tau_{2^{s+1},2^{s+1}}\bigr)
\]
and $M_s$ be the Gram matrix of the ridge functions associated with $\Theta_{2^{s+1}}$. Using $h_{2^{s+1}}\simeq J_s^{-1}$, the one-sided estimate proved in \cite[(4.34)]{liu2025condition} gives
\begin{equation}\label{eqn:item3_endpoint_weighted_bound}
\bigl(v^{(s)}\bigr)^\top W_s^{1/2}M_sW_s^{1/2}v^{(s)}
\lesssim J_s^{-(d+2k+1)}.
\end{equation}

Set
\[
a_i^{(s)}:=\sqrt{\tau_{i,2^{s+1}}}\,v_i^{(s)},\qquad 2^s<i\leq2^{s+1},
\]
then $\tau_{i,2^{s+1}}\simeq2^{-(s+1)}$ and $\|v^{(s)}\|_2=1$ yields
\begin{equation}\label{eqn:item3_endpoint_beta}
\|a_i^{(s)}\|_2^2\simeq2^{-(s+1)}.
\end{equation}
Denote
\[
g_s:=2^{\frac{s+1}{2}}\sum_{i=2^s+1}^{2^{s+1}}a_i^{(s)}\sigma_k(\theta_i^*\cdot\circ).
\]
Then for $m\in E_{\sigma_k}$, $m\leq J_s$ and $1\leq\ell\leq N(m)$, \eqref{eqn:fn_harmonic_coeff} and $v^{(s)}\in\mathcal V(J_s)$ give
\[
\begin{split}
\widehat{g_s}(m,\ell)
&=2^{\frac{s+1}{2}}\widehat{\sigma_k}(m)\sum_{i=2^s+1}^{2^{s+1}}a_i^{(s)}Y_{m,\ell}(\theta_i^*)=2^{\frac{s+1}{2}}\widehat{\sigma_k}(m)\bigl(v^{(s)}\bigr)^\top V_{2^{s+1}}(m,\ell)=0.
\end{split}
\]
Recall $\widehat{\sigma_k}(m)$ is supported in $E_{\sigma_k}$,
\begin{equation}\label{eqn:item3_endpoint_low_frequency}
\mathcal P_{J_s}g_s=0.
\end{equation}
Moreover, \eqref{eqn:item3_endpoint_weighted_bound} yields
\begin{equation}\label{eqn:item3_endpoint_norm}
\begin{split}
\|g_s\|_2^2
&=2^{s+1}\bigl(v^{(s)}\bigr)^\top W_s^{1/2}M_sW_s^{1/2}v^{(s)}\lesssim2^{s+1}J_s^{-(d+2k+1)}
\lesssim J_s^{-(2k+1)}.
\end{split}
\end{equation}

Define
\[
f:=\sum_{s=s_0}^{\infty}g_s,\qquad F_s:=\sum_{t=s_0}^sg_t.
\]
By \eqref{eqn:item3_endpoint_norm}, then $f\in\mathcal L^2(\SS^d)$ and by nestedness, $F_s\in L_{2^{s+1}}^k(\Theta_{2^{s+1}})$. Now for $2^{s+1}\leq n<2^{s+2}$, nestedness and \eqref{eqn:item3_endpoint_norm} yield \eqref{eqn:rate_alpha}:
\[
\begin{split}
\inf_{v_n\in L_n^k(\Theta_n)}\|f-v_n\|_2
&\leq\|f-F_s\|_2
\leq\sum_{t>s}\|g_t\|_2\lesssim\sum_{t>s}2^{-\frac{(t+1)}{d}\frac{2k+1}{2}}
\lesssim2^{-\frac{(s+2)}{d}\frac{2k+1}{2}}
\lesssim n^{-\frac{2k+1}{2d}}.
\end{split}
\]

It remains to prove that $f\notin\mathcal B_{2,\infty}^{r}(\SS^d)$. When $F_s$ is represented over $\Theta_{2^{s+1}}$, the disjointness of the coefficient supports and \eqref{eqn:item3_endpoint_beta} imply
\[
\sum_{t=s_0}^s\sum_{i=2^t+1}^{2^{t+1}}2^{t+1}\bigl|a_i^{(t)}\bigr|^2\simeq s-s_0+1\simeq s.
\]
Let $\displaystyle q_s=\Big\lceil\log_2\Big(\frac{C_3}{\underline h_{2^{s+1}}}\Big)\Big\rceil$, applying \cite[Lemma~2.1]{mao2025sharp} therefore gives
\begin{equation}\label{eqn:item3_endpoint_partial_tail}
\|(\mathcal P_{2^{q_s+2}}-\mathcal P_{2^{q_s}})F_s\|_2
\gtrsim2^{-q_s(k+\frac12)}\sqrt{s}.
\end{equation}
On the other hand, antipodally quasi-uniformity yields $2^{q_s}\simeq\underline{h}^{-1}\simeq2^{\frac{s}{d}}$, \eqref{eqn:item3_endpoint_norm} implies
\[
\|(\mathcal P_{2^{q_s+2}}-\mathcal P_{2^{q_s}})(f-F_s)\|_2
\leq\|f-F_s\|_2
\lesssim2^{-\frac{s+2}{d}(k+\frac12)}
\lesssim2^{-q_s(k+\frac12)}.
\]
Combining this with \eqref{eqn:item3_endpoint_partial_tail}, we obtain
\[
2^{q_s(k+\frac12)}\|(\mathcal P_{2^{q_s+2}}-\mathcal P_{2^{q_s}})f\|_2
\gtrsim \sqrt{s}\longrightarrow\infty,
\]
which implies $f\notin\mathcal B_{2,\infty}^{r}(\SS^d)$.

Finally, every $g_s$ has the required parity above degree $k$, and this property is preserved under convergence in $\mathcal L^2(\SS^d)$. Therefore
\[
f\in\mathcal L_*^2(\SS^d)\setminus\mathcal B_{2,\infty,*}^{k+\frac{1}{2}}(\SS^d).
\]

\subsection{Proof of item \ref{item4}}
By \eqref{eqn:item4_two_rates}, for every $n$ we can choose
\[
    f_n=\sum_{j=1}^na_{j,n}\sigma_k(\theta_{j,n}^*\cdot\circ),\qquad \widetilde f_n=\sum_{j=1}^n\widetilde a_{j,n}\sigma_k(\widetilde\theta_{j,n}^*\cdot\circ)
\]
such that
\begin{equation}\label{eqn:item4_approx}
    \|f-f_n\|_{\mathcal L^2(\SS^d)}+\|f-\widetilde f_n\|_{\mathcal L^2(\SS^d)}\lesssim n^{-\frac{d+2k+1}{2d}+\alpha}.
\end{equation}
Applying Theorem~\ref{lem:ridge_stability}, \eqref{eqn:item4_union_separation}, and the coefficients $a_{1,n},\ldots,a_{n,n},-\widetilde a_{1,n},\ldots,-\widetilde a_{n,n}$ gives
\begin{equation}\label{eqn:item4_coefficient_lower}
    \|f_n-\widetilde f_n\|_{\mathcal L^2(\SS^d)}^2\gtrsim n^{-\frac{2k+1}{d}}\Big(\sum_{j=1}^na_{j,n}^2+\sum_{j=1}^n\widetilde a_{j,n}^2\Big).
\end{equation}
On the other hand, \eqref{eqn:item4_approx} yields
\[
    \|f_n-\widetilde f_n\|_{\mathcal L^2(\SS^d)}\leq\|f-f_n\|_{\mathcal L^2(\SS^d)}+\|f-\widetilde f_n\|_{\mathcal L^2(\SS^d)}\lesssim n^{-\frac{d+2k+1}{2d}+\alpha}.
\]
Thus we obtain
\begin{equation}\label{eqn:item4_coefficient_bound}
    \sum_{j=1}^na_{j,n}^2+\sum_{j=1}^n\widetilde a_{j,n}^2\lesssim n^{2\alpha-1}.
\end{equation}
and in particular,
\[
    \Big(n\sum_{j=1}^na_{j,n}^2\Big)^{\frac{1}{2}}\lesssim n^\alpha.
\]
Therefore, for some $M>0$ independent of $n$, $f_n\in L_{n,M}^{k,\alpha}(\Theta_n)$, while \eqref{eqn:item4_approx} retains the rate in item~\ref{item1}. The inverse implication in item~\ref{item1} therefore gives $f\in\mathcal B_{2,\infty,*}^{\frac{d+2k+1}{2}-\alpha d}(\SS^d)$, proving item~\ref{item4}.

\bibliographystyle{abbrv}
\bibliography{ref}

@article{liu2025condition,
  title={Condition numbers and eigenvalue spectra of shallow networks on spheres},
  author={Liu, Xinliang and Mao, Tong and Xu, Jinchao},
  journal={arXiv preprint arXiv:2511.02625},
  year={2025}
}

@article{E2019,
  title={A priori estimates and analysis of a numerical method for machine learning},
  author={E, Weinan and Yu, Bing},
  journal={Science China Mathematics},
  volume={62},
  number={11},
  pages={2303--2330},
  year={2019},
  publisher={Springer}
}

@book{dai2013approximation,
  title={Approximation theory and harmonic analysis on spheres and balls},
  author={Dai, Feng and Xu, Yuan},
  year={2013},
  publisher={Springer}
}

@article{mhaskar2020kernel,
  title={Kernel-based analysis of massive data},
  author={Mhaskar, Hrushikesh N},
  journal={Frontiers in Applied Mathematics and Statistics},
  volume={6},
  pages={30},
  year={2020},
  publisher={Frontiers Media SA}
}

@article{zhou2006approximation,
  title={Approximation with polynomial kernels and SVM classifiers},
  author={Zhou, Ding-Xuan and Jetter, Kurt},
  journal={Advances in Computational Mathematics},
  volume={25},
  number={1},
  pages={323--344},
  year={2006},
  publisher={Springer}
}

@article{maiorov2006pseudo,
  title={Pseudo-dimension and entropy of manifolds formed by affine-invariant dictionary},
  author={Maiorov, Vitaly},
  journal={Advances in Computational Mathematics},
  volume={25},
  pages={435--450},
  year={2006},
  publisher={Springer}
}

@article{devore1989optimal,
  title={Optimal nonlinear approximation},
  author={DeVore, Ronald A, Ralph Howard and C. Micchelli},
  journal={Manuscripta mathematica},
  volume={58},
  number={3},
  pages={295--306},
  year={1989},
  publisher={Springer}
}

@article{guo2020optimal,
  title={Optimal approximation with ReLU neural networks},
  author={Guo, Zheng-Chu and Zhou, Ding-Xuan},
  journal={Neural Networks},
  volume={128},
  pages={96--109},
  year={2020},
  publisher={Elsevier}
}

@article{bernstein2020distance,
  title={On the distance between two neural networks and the stability of learning},
  author={Bernstein, Jeremy and Vahdat, Arash and Yue, Yisong and Liu, Ming-Yu},
  journal={Advances in Neural Information Processing Systems},
  volume={33},
  pages={21370--21381},
  year={2020}
}

@article{huang2023generalization,
  title={Generalization analysis of pairwise learning for ranking with deep neural networks},
  author={Huang, Shuo and Zhou, Junyu and Feng, Han and Zhou, Ding-Xuan},
  journal={Neural Computation},
  volume={35},
  number={6},
  pages={1135--1158},
  year={2023},
  publisher={MIT Press One Rogers Street, Cambridge, MA 02142-1209, USA journals-info~…}
}

@article{bach2017breaking,
  title={Breaking the curse of dimensionality with convex neural networks},
  author={Bach, Francis},
  journal={The Journal of Machine Learning Research},
  volume={18},
  number={1},
  pages={629--681},
  year={2017},
  publisher={JMLR. org}
}

@article{mhaskar2006weighted,
  title={Weighted quadrature formulas and approximation by zonal function networks on the sphere},
  author={Mhaskar, Hrushikesh Narhar},
  journal={Journal of Complexity},
  volume={22},
  number={3},
  pages={348--370},
  year={2006},
  publisher={Elsevier}
}

@article{schneider1967problem,
  title={Zu einem Problem von Shephard {\"u}ber die Projektionen konvexer K{\"o}rper},
  author={Schneider, Rolf},
  journal={Mathematische Zeitschrift},
  volume={101},
  pages={71--82},
  year={1967},
  publisher={Springer}
}

@article{narcowich2006sobolev,
  title={Sobolev error estimates and a Bernstein inequality for scattered data interpolation via radial basis functions},
  author={Narcowich, Francis J and Ward, Joseph D and Wendland, Holger},
  journal={Constructive Approximation},
  volume={24},
  pages={175--186},
  year={2006},
  publisher={Springer}
}

@article{mhaskar2001spherical,
  title={Spherical Marcinkiewicz-Zygmund inequalities and positive quadrature},
  author={Mhaskar, H and Narcowich, F and Ward, J},
  journal={Mathematics of computation},
  volume={70},
  number={235},
  pages={1113--1130},
  year={2001}
}

@book{devore1993constructive,
  title={Constructive approximation},
  author={DeVore, Ronald A and Lorentz, George G},
  volume={303},
  year={1993},
  publisher={Springer Science \& Business Media}
}

@article{siegel2022sharp,
  title={Sharp bounds on the approximation rates, metric entropy, and n-widths of shallow neural networks},
  author={Siegel, Jonathan W and Xu, Jinchao},
  journal={Foundations of Computational Mathematics},
  pages={1--57},
  year={2022},
  publisher={Springer}
}

@article{siegel2023characterization,
  title={Characterization of the variation spaces corresponding to shallow neural networks},
  author={Siegel, Jonathan W and Xu, Jinchao},
  journal={Constructive Approximation},
  volume={57},
  number={3},
  pages={1109--1132},
  year={2023},
  publisher={Springer}
}

@article{mhaskar1993approximation,
  title={Approximation properties of a multilayered feedforward artificial neural network},
  author={Mhaskar, Hrushikesh Narhar},
  journal={Advances in Computational Mathematics},
  volume={1},
  pages={61--80},
  year={1993},
  publisher={Springer}
}

@article{pinkus1999approximation,
  title={Approximation theory of the MLP model in neural networks},
  author={Pinkus, Allan},
  journal={Acta numerica},
  volume={8},
  pages={143--195},
  year={1999},
  publisher={Cambridge University Press}
}

@article{fang2020theory,
  title={Theory of deep convolutional neural networks II: Spherical analysis},
  author={Fang, Zhiying and Feng, Han and Huang, Shuo and Zhou, Ding-Xuan},
  journal={Neural Networks},
  volume={131},
  pages={154--162},
  year={2020},
  publisher={Elsevier}
}

@article{devore1996some,
  title={Some remarks on greedy algorithms},
  author={DeVore, Ronald A and Temlyakov, Vladimir N},
  journal={Advances in computational Mathematics},
  volume={5},
  number={1},
  pages={173--187},
  year={1996},
  publisher={Springer}
}

@article{siegel2022optimal,
  title={Optimal convergence rates for the orthogonal greedy algorithm},
  author={Siegel, Jonathan W and Xu, Jinchao},
  journal={IEEE Transactions on Information Theory},
  volume={68},
  number={5},
  pages={3354--3361},
  year={2022},
  publisher={IEEE}
}

@article{petrushev1998approximation,
  title={Approximation by ridge functions and neural networks},
  author={Petrushev, Pencho P},
  journal={SIAM Journal on Mathematical Analysis},
  volume={30},
  number={1},
  pages={155--189},
  year={1998},
  publisher={SIAM}
}

@article{klusowski2018approximation,
  title={Approximation by combinations of ReLU and squared ReLU ridge functions with l1 and l0 controls},
  author={Klusowski, Jason M and Barron, Andrew R},
  journal={IEEE Transactions on Information Theory},
  volume={64},
  number={12},
  pages={7649--7656},
  year={2018},
  publisher={IEEE},
  doi={10.1109/TIT.2018.2834518}
}

@article{barron1993universal,
  title={Universal approximation bounds for superpositions of a sigmoidal function},
  author={Barron, Andrew R},
  journal={IEEE Transactions on Information theory},
  volume={39},
  number={3},
  pages={930--945},
  year={1993},
  publisher={IEEE}
}

@inproceedings{bourgain2006projection,
  title={Projection bodies},
  author={Bourgain, Jean and Lindenstrauss, Joram},
  booktitle={Geometric Aspects of Functional Analysis: Israel Seminar (GAFA) 1986--87},
  pages={250--270},
  year={2006},
  organization={Springer}
}

@article{barron1994approximation,
  title={Approximation and estimation bounds for artificial neural networks},
  author={Barron, Andrew R},
  journal={Machine learning},
  volume={14},
  pages={115--133},
  year={1994},
  publisher={Springer}
}

@article{siegel2022high,
  title={High-order approximation rates for shallow neural networks with cosine and {ReLU}k activation functions},
  author={Siegel, Jonathan W and Xu, Jinchao},
  journal={Applied and Computational Harmonic Analysis},
  volume={58},
  pages={1--26},
  year={2022},
  publisher={Elsevier}
}

@article{mhaskar2023tractability,
  title={Tractability of approximation by general shallow networks},
  author={Mhaskar, Hrushikesh and Mao, Tong},
  journal={arXiv preprint arXiv:2308.03230},
  year={2023}
}

@article{liu2025integral,
  title={Integral Representations of Sobolev Spaces via {ReLU}$^k$ Activation Function and Optimal Error Estimates for Linearized Networks},
  author={Liu, Xinliang and Mao, Tong and Xu, Jinchao},
  journal={arXiv preprint arXiv:2505.00351},
  year={2025}
}

@article{Jones1992,
  author = {Lee K Jones},
  title = {A Simple Lemma on Greedy Approximation in Hilbert Space and Convergence Rates for Projection Pursuit Regression and Neural Network Training},
  journal = {The Annals of Statistics},
  volume = {20},
  number = {1},
  pages = {608-613},
  year = {1992}
}

@article{bernstein1912quelques,
  author    = {S. N. Bernstein},
  title     = {Sur l'ordre de la meilleure approximation des fonctions continues par les polyn{\^o}mes de degr{\'e} donn{\'e}.},
  journal   = {Mem. Cl. Sci. Acad. Roy. Belg.},
  year      = {1912},
  volume    = {4},
  pages     = {1--103}
}

@book{nikolskii1975approximation,
  author    = {S. M. Nikolskii},
  title     = {Approximation of Functions of Several Variables and Imbedding Theorems},
  publisher = {Springer},
  year      = {1975},
  address   = {Berlin}
}

@article{mhaskar2004polynomialLocal,
  author  = {Mhaskar, H. N.},
  title   = {Polynomial operators and local smoothness classes on the unit interval},
  journal = {Journal of Approximation Theory},
  volume  = {131},
  number  = {2},
  pages   = {243--267},
  year    = {2004}
}

@article{mhaskar2005sphereLocalized,
  author  = {Mhaskar, H. N.},
  title   = {On the representation of smooth functions on the sphere using finitely many bits},
  journal = {Applied and Computational Harmonic Analysis},
  volume  = {18},
  number  = {3},
  pages   = {215--233},
  year    = {2005}
}

@article{mhaskar2019reluZonal,
  author  = {Mhaskar, H. N.},
  title   = {Function approximation with zonal function networks with activation functions analogous to the rectified linear unit functions},
  journal = {Journal of Complexity},
  volume  = {51},
  pages   = {1--19},
  year    = {2019}
}

@inproceedings{cohen2018spherical,
  title={Spherical {CNNs}},
  author={Cohen, Taco S. and Geiger, Mario and K{\"o}hler, Jonas and Welling, Max},
  booktitle={International Conference on Learning Representations},
  year={2018}
}

@inproceedings{jiang2019spherical,
  title={Spherical {CNNs} on Unstructured Grids},
  author={Jiang, Chiyu and Huang, Jingwei and Kashinath, Karthik and Prabhat and Marcus, Philip and Niessner, Matthias},
  booktitle={International Conference on Learning Representations},
  year={2019}
}

@article{mao2025sharp,
  title={Sharp Lower Bounds for Linearized ReLU\^{} k Approximation on the Sphere},
  author={Mao, Tong and Xu, Jinchao},
  journal={arXiv preprint arXiv:2510.04060},
  year={2025}
}

@inproceedings{defferrard2019deepsphere,
  title={{DeepSphere}: Towards an Equivariant Graph-Based Spherical {CNN}},
  author={Defferrard, Micha{\"e}l and Perraudin, Nathana{\"e}l and Kacprzak, Tomasz and Sgier, Raphael},
  booktitle={ICLR Workshop on Representation Learning on Graphs and Manifolds},
  year={2019},
  eprint={1904.05146},
  archivePrefix={arXiv}
}

@article{feng2023generalization,
  title={Generalization Analysis of {CNNs} for Classification on Spheres},
  author={Feng, Han and Huang, Shuo and Zhou, Ding-Xuan},
  journal={IEEE Transactions on Neural Networks and Learning Systems},
  year={2023}
}

@article{yang2026spherical,
  title={Spherical Analysis of Learning Nonlinear Functionals},
  author={Yang, Yang and Huang, Shuo and Feng, Han and Zhou, Ding-Xuan},
  journal={Journal of Machine Learning Research},
  year={2026}
}

@article{liu2026optimal,
  title={Optimal Neural Network Approximation via Empirical Least Squares with Deterministic Samples},
  author={Liu, Xinliang and Mao, Tong and Xu, Jinchao},
  journal={arXiv preprint arXiv:2608.06687},
  year={2026}
}

@article{devore1998nonlinear,
  author  = {DeVore, Ronald A.},
  title   = {Nonlinear Approximation},
  journal = {Acta Numerica},
  volume  = {7},
  pages   = {51--150},
  year    = {1998},
  doi     = {10.1017/S0962492900002816}
}

\end{document}